\documentclass[a4paper]{article}
\usepackage{hyperref}
\usepackage{amsmath}
\usepackage{amsfonts}
\usepackage{amsthm}
\usepackage{indentfirst}
\usepackage{geometry}

\newtheorem{theorem}{Theorem}

\catcode`@=11
\@addtoreset{equation}{section}
\catcode`@=12

\begin{document}

\begin{center}
{\Large  A Liouville theorem for the fractional Ginzburg-Landau equation}
\end{center}

\vskip 5mm

\begin{center}
{\sc Yayun Li} \\
\vskip 3mm
School of Applied Mathematics\\
Nanjing University of Finance \& Economics,
Nanjing, 210023, China\\
Email:yayunli@nufe.edu.cn
\end{center}

\begin{center}
{\sc Qinghua Chen and Yutian Lei} \\
\vskip 3mm
Institute of Mathematics\\
School of Mathematical Sciences\\
Nanjing Normal University,
Nanjing, 210023, China\\
Email:leiyutian@njnu.edu.cn
\end{center}

\vskip 5mm {\leftskip5mm\rightskip5mm \normalsize
\noindent{\bf{Abstract}}
In this paper, we are concerned with a Liouville-type
result of the nonlinear integral equation
\begin{equation*}
u(x)=\int_{\mathbb{R}^{n}}\frac{u(1-|u|^{2})}{|x-y|^{n-\alpha}}dy,
\end{equation*}
where $u: \mathbb{R}^{n} \to \mathbb{R}^{k}$ with $k \geq 1$ and $1<\alpha<n/2$.
We prove that $u \in L^2(\mathbb{R}^n) \Rightarrow u \equiv 0$ on $\mathbb{R}^n$,
as long as $u$ is a bounded and differentiable solution.

\par
\noindent{\bf{Keywords}}: Ginzburg-Landau equation, Liouville theorem, fractional
Laplacian
\par
{\bf{MSC2010}}: 45G05, 45E10, 35Q56, 35R11}

\vskip 1cm

If a harmonic function $u$ is bounded on $\mathbb{R}^n$, then $u \equiv Const.$
(this is the Liouville theorem).
Moreover, if $u$ is integrable
(i.e. $u \in L^s(\mathbb{R}^n)$ for some $s \geq 1$), then $u \equiv 0$ on $\mathbb{R}^n$.

In 1994, Brezis, Merle and Rivi\`{e}re \cite{BMR} studied the
quantization effects of the following equation
\begin{equation}\label{1.5}
-\Delta u=u(1-|u|^2) \quad on \ \mathbb{R}^2.
\end{equation}
Here $u:\mathbb{R}^2 \to \mathbb{R}^2$ is a vector valued function.
It is the Euler-Lagrange equation of the Ginzburg-Landau energy
$$
E_{GL}(u)=\frac{1}{2}\|\nabla
u\|_{L^2(\mathbb{R}^2)}^2+\frac{1}{4}\|1-|u|^2\|_{L^2(\mathbb{R}^2)}^2.
$$
In particular,
they proved the finite energy solution (i.e., $u$ satisfies $\nabla
u \in L^2(\mathbb{R}^2)$) is bounded (see also \cite{HH} and \cite{Ma2})
\begin{equation}\label{1.7}
|u| \leq 1 \quad on \ \mathbb{R}^n.
\end{equation}
(Here $n=2$.) Based on this result, they obtained a Liouville type theorem for finite
energy solutions (cf. Theorem 2 in \cite{BMR}):

\emph{Let $u:\mathbb{R}^2 \to \mathbb{R}^2$ be a classical solution of (\ref{1.5}).
If $\nabla u \in L^2(\mathbb{R}^2)$, then either $u \in L^2(\mathbb{R}^2)$ which implies
$u \equiv 0$, or $1-|u|^2 \in L^1(\mathbb{R}^2)$ which implies $u \equiv C$ with
$|C|=1$.}

The boundedness and the integrability of solutions are the important conditions
which ensure that the Liouville theorem holds.
The Pohozaev identity plays a key role in the proof.

In this paper, we are concerned with the integral equation
\begin{equation}\label{1.1}
u(x)=\int_{\mathbb{R}^n}\frac{u(1-|u|^2)}{|x-y|^{n-\alpha}}dy.
\end{equation}
Here $u:\mathbb{R}^n \to \mathbb{R}^k$, $k \geq 1$, $n \geq 3$, and $1<\alpha<n/2$.
We also apply the integral form of the Pohozaev identity (which
was used for the Lane-Emden equations in \cite{CDM}, \cite{LeiL}
and \cite{Xu}) to establish a Liouville theorem.

\begin{theorem} \label{th1.1}
Assume that $u:\mathbb{R}^n \to \mathbb{R}^k$ is bounded and differentiable, and solves (\ref{1.1})
with $\alpha \in (1,n/2)$. If $u \in L^2(\mathbb{R}^n)$, then $u(x) \equiv 0$.
\end{theorem}

\begin{proof}
For convenience, we denote $B_R(0)$ by $B_R$ here.

{\it Step 1.} We claim that the  improper integral
\begin{equation}\label{jia}
\int_{\mathbb{R}^{n}}\frac{z\cdot\nabla[u(z)(1-|u(z)|^{2})]}{|x-z|^{n-\alpha}}dz
\end{equation}
is convergent at each $x \in \mathbb{R}^n$.

In fact, since $u \in L^2(\mathbb{R}^n)$, we can find $R=R_j \to \infty$
such that
\begin{equation}\label{yi}
R\int_{\partial B_{R}}|u(z)|^{2}ds \to 0.
\end{equation}
Since $u$ is bounded, by the H\"older inequality, we obtain that for
sufficiently large $R$, there holds
$$\begin{array}{ll}
\displaystyle R\left|\int_{\partial B_R} \frac{u(z)(1-|u(z)|^{2})}{|x-z|^{n-\alpha}}ds \right|
&\displaystyle \leq CR^{1-n+\alpha} \int_{\partial B_R}|u(z)|ds\\[3mm]
&\displaystyle \leq CR^{1-n+\alpha} \left(R\int_{\partial B_{R}}|u(z)|^{2}ds\right)^{\frac{1}{2}}
R^{-\frac{1}{2}+\frac{n-1}{2}}.
\end{array}
$$
Let $R=R_j \to \infty$. Noting $\alpha < n/2$, and using (\ref{yi}) we get
\begin{equation}\label{bing}
R\int_{\partial B_R} \frac{u(z)(1-|u(z)|^{2})}{|x-z|^{n-\alpha}}ds
\to 0
\end{equation}
when $R=R_j \to \infty$.

Next, we claim that the improper integral
\begin{equation}\label{ji}
I(\mathbb{R}^n):=\int_{\mathbb{R}^n} \frac{u(z)(1-|u(z)|^{2})(x-z)\cdot z}
{|x-z|^{n-\alpha+2}}dz
\end{equation}
absolutely converges for each $x \in \mathbb{R}^n$.

In fact, we observe that the defect points of $I(\mathbb{R}^n)$ are $x$ and $\infty$.
When $z$ is near $\infty$, we have
$$
|I(\mathbb{R}^n\setminus B_r)| \leq C\int_{\mathbb{R}^n\setminus B_r}\frac{|u(z)|dz}{|x-z|^{n-\alpha}}
\leq C\left(\int_{\mathbb{R}^n}|u|^2dz\right)^{\frac{1}{2}}
\left(\int_r^\infty \rho^{n-2(n-\alpha)} \frac{d\rho}{\rho}\right)^{\frac{1}{2}}.
$$
In view of $u \in L^2(\mathbb{R}^n)$ and $\alpha<n/2$, we get
\begin{equation}\label{xin}
|I(\mathbb{R}^n\setminus B_r)|<\infty.
\end{equation}
When $z$ is near $x$, we first take
\begin{equation}\label{ding}
s \in \left(\frac{n}{\alpha-1},\infty\right).
\end{equation}
Clearly, $1<\alpha<n/2$ implies $s>2$. In addition,
\begin{equation}\label{wu}
u \in L^s(\mathbb{R}^n)
\end{equation}
because $u$ is bounded and $u \in L^2(\mathbb{R}^n)$.
Note that
$$
|I(B_\delta(x))| \leq C\int_{B_\delta(x)}\frac{|u(z)|dz}{|x-z|^{n-\alpha+1}}
\leq C\left(\int_{\mathbb{R}^n}|u|^sdz\right)^{\frac{1}{s}}
\left(\int_0^r \rho^{n-\frac{s}{s-1}(n-\alpha+1)} \frac{d\rho}{\rho}\right)^{1-\frac{1}{s}}.
$$
By (\ref{ding}) and (\ref{wu}), we get
$$
|I(B_\delta(x))|<\infty.
$$
Combining this with (\ref{xin}), we prove that (\ref{ji}) is absolutely convergent.

Finally we prove that (\ref{jia}) is convergent. Integrating by parts yields
\begin{equation}\label{geng}
\begin{array}{ll}
&\displaystyle\int_{B_R}\frac{z\cdot\nabla[u(z)(1-|u(z)|^{2})]}{|x-z|^{n-\alpha}}dz
=R\int_{\partial B_R} \frac{u(z)(1-|u(z)|^{2})}{|x-z|^{n-\alpha}}ds\\[3mm]
&\quad \displaystyle -n\int_{B_R}\frac{u(z)(1-|u(z)|^2)}{|x-z|^{n-\alpha}}dz
-(n-\alpha)\int_{B_R} \frac{u(z)(1-|u(z)|^{2})(x-z)\cdot z}{|x-z|^{n-\alpha+2}}dz.
\end{array}
\end{equation}
Letting $R=R_j \to \infty$ in (\ref{geng}) and using (\ref{1.1}) and (\ref{bing}),
we can see that
$$
\int_{\mathbb{R}^n}\frac{z\cdot\nabla[u(z)(1-|u(z)|^{2})]}{|x-z|^{n-\alpha}}dz
=-nu(x)+(\alpha-n)I(\mathbb{R}^n),
$$
and hence it is convergent at each $x \in \mathbb{R}^n$.

{\it Step 2.} Proof of Theorem \ref{th1.1}.

For any $\lambda>0$, from (\ref{1.1}) it follows
$$
u(\lambda x)=\lambda^{\alpha}\int_{\mathbb{R}^{n}}\frac{u(\lambda z)
(1-|u(\lambda z)|^{2})}{|x-z|^{n-\alpha}}dz.
$$
Differentiating both sides with respect to $\lambda$ yields
$$\begin{array}{ll}
x\cdot\nabla u(\lambda x)
&=\displaystyle\alpha\lambda^{\alpha-1}\int_{\mathbb{R}^{n}}\frac{u(\lambda z)
(1-|u(\lambda z)|^{2})}{|x-z|^{n-\alpha}}dz\\[3mm]
&\quad +\displaystyle\lambda^{\alpha}\int_{\mathbb{R}^{n}}\frac{(z\cdot \nabla u(\lambda z))
(1-|u(\lambda z)|^{2})+u(\lambda z)
[-2u(\lambda z)(z\cdot \nabla  u(\lambda z))]}{|x-z|^{n-\alpha}}dz.
\end{array}
$$
Letting $\lambda=1$ yields
\begin{equation}\label{d}
x\cdot\nabla u(x)
=\alpha u(x)+\int_{\mathbb{R}^{n}}\frac{z\cdot\nabla[u(1-|u|^{2})]}{|x-z|^{n-\alpha}}dz.
\end{equation}

Since $u$ is bounded and $u\in L^2(\mathbb{R}^n)$, it follows that $u \in L^4(\mathbb{R}^n)$, and hence
\begin{equation}\label{ren}
R\int_{\partial B_{R}}|u|^{4}ds \to 0
\end{equation}
for some $R=R_j \to \infty$. Thus, integrating by parts and using (\ref{yi})
and (\ref{ren}), we respectively obtain
\begin{equation}\label{e}
\int_{\mathbb{R}^{n}}u(x)(x\cdot\nabla u(x))dx=\frac{-n}{2}\int_{\mathbb{R}^{n}}|u(x)|^{2}dx,
\end{equation}
and
\begin{equation}\label{f}
\int_{\mathbb{R}^{n}}u(x)|u(x)|^{2}(x\cdot\nabla u(x))dx=\frac{-n}{4}\int_{\mathbb{R}^{n}}|u(x)|^{4}dx.
\end{equation}
These results show that
\begin{equation}\label{gui}
\int_{\mathbb{R}^{n}}u(x)(1-|u(x)|^{2})(x\cdot\nabla u(x))dx<\infty.
\end{equation}

Multiply (\ref{d}) by $u(x)(1-|u(x)|^{2})$ and integrate over $B_R$.
Letting $R=R_j \to \infty$, from $u \in L^2(\mathbb{R}^n) \cap L^4(\mathbb{R}^n)$ and
(\ref{gui}), we get
$$
\int_{\mathbb{R}^{n}}u(x)(1-|u(x)|^{2})
\int_{\mathbb{R}^{n}}\frac{z\cdot\nabla[u(1-|u|^{2})]}{|x-z|^{n-\alpha}}dzdx
<\infty,
$$
and
\begin{equation}\label{g}
\begin{array}{ll}
&\displaystyle\int_{\mathbb{R}^{n}}u(x)(1-|u(x)|^{2})(x\cdot\nabla u(x))dx
-\displaystyle\alpha\int_{\mathbb{R}^{n}}|u(x)|^{2}(1-|u(x)|^{2})dx\\[3mm]
&=\displaystyle\int_{\mathbb{R}^{n}}u(x)(1-|u(x)|^{2})
\int_{\mathbb{R}^{n}}\frac{z\cdot\nabla[u(1-|u|^{2})]}{|x-z|^{n-\alpha}}dzdx.
\end{array}
\end{equation}

We use the Fubini theorem and (\ref{1.1}) to handle the right hand side term.
Thus,
\begin{equation}\label{g1}
\begin{array}{ll}
&\quad \displaystyle\int_{\mathbb{R}^{n}}u(x)(1-|u(x)|^{2})\int_{\mathbb{R}^{n}}
\frac{z\cdot\nabla[u(1-|u|^{2})]}{|x-z|^{n-\alpha}}dzdx\\[3mm]
&=\displaystyle\int_{\mathbb{R}^{n}}z\cdot\nabla[u(1-|u|^{2})]
\int_{\mathbb{R}^{n}}\frac{u(x)(1-|u(x)|^{2})}{|x-z|^{n-\alpha}}dxdz\\[3mm]
&=\displaystyle\int_{\mathbb{R}^{n}}(x\cdot\nabla[u(1-|u|^{2})])u(x)dx\\[3mm]
&=\displaystyle\int_{\mathbb{R}^{n}}u(x)(1-|u|^{2})(x\cdot\nabla u(x))dx
-\int_{\mathbb{R}^{n}}|u|^{2}(x\cdot\nabla |u(x)|^2)dx.
\end{array}
\end{equation}

Inserting this result into (\ref{g}), and using (\ref{f}) we have
\begin{equation}\label{h}
\int_{\mathbb{R}^{n}}[\alpha |u(x)|^{2}+(\frac{n}{2}-\alpha)|u(x)|^{4}]dx=0.
\end{equation}
In view of $\alpha \in (1, n/2)$, (\ref{h}) leads to $|u|\equiv 0.$
Theorem \ref{th1.1} is proved.
\end{proof}

\paragraph{Remark 1.}
Clearly, (\ref{h}) implies
\begin{equation}\label{i}
\int_{\mathbb{R}^{n}}[\alpha (|u(x)|^{2}-|u(x)|^{4})+\frac{n}{2}|u(x)|^{4}]dx=0.
\end{equation}
When $u$ satisfies (\ref{1.7}), (\ref{i}) also implies $|u|\equiv 0.$
For (\ref{1.5}),  the bound $|u| \leq 1$ for solutions
$u : \mathbb{R}^{n} \to \mathbb{R}^{k}$ was first proved by Brezis (cf. \cite{Brezis}).
Ma also pointed out that (\ref{1.7}) holds true (cf. \cite{Ma2}).

\paragraph{Remark 2.}
In 2016, Ma \cite{Ma1} proved (\ref{1.7}) for the
Ginzburg-Landau-type equation with fractional Laplacian
\begin{equation}\label{1.2}
(-\Delta)^{\frac{\alpha}{2}} u=(1-|u|^{2})u  \quad on \ \mathbb{R}^n
\end{equation}
under the assumption
\begin{equation}\label{zi}
1-|u|^2 \in L^2(\mathbb{R}^n),
\end{equation}
where $n \geq 2$ and $0<\alpha<2$.
The physical background of (\ref{1.2}) can be found in \cite{MR} and \cite{TZ}.
Such an equation with $\alpha=1$ was well studied in \cite{MS}.
Recall the definition of fractional Laplacian on $\mathbb{R}^{n}$. Let
$n \geq 2$ and $0<\alpha<2$. Write
$$
E=L_{\alpha} \cap C_{loc}^{1,1}(\mathbb{R}^{n}),
$$
where $L_{\alpha}=\left\{u\in L_{loc}^{1}(\mathbb{R}^{n});
\int_{\mathbb{R}^{n}}\frac{|u(x)|dx}{1+|x|^{n+\alpha}}<\infty\right\}$. For
a vector value function $u \in E$ from $\mathbb{R}^n$ to $\mathbb{R}^k$, define
\begin{equation}\label{qz}
(-\triangle)^{\frac{\alpha}{2}}u:=C_{n,\alpha}P.V. \int_{\mathbb{R}^n}
\frac{u(x)-u(y)}{|x-y|^{n+\alpha}}dy
=C_{n,\alpha}\lim_{\varepsilon \to 0^+} \int_{|x-y| \geq \varepsilon}
\frac{u(x)-u(y)}{|x-y|^{n+\alpha}}dy.
\end{equation}
Here $C_{n,\alpha}$ is a positive constant.

Clearly, (\ref{zi}) and $u \in L^2(\mathbb{R}^n)$ are incompatible.

\paragraph{Remark 3.}
Another definition of the fractional order Laplacian involves
the Riesz potential (cf. Chapter 5 in \cite{Stein}).
Assume $\alpha \in (0,n)$, and $u, u(1-|u|^2) \in \mathcal{S}'(\mathbb{R}^n)$,
then (\ref{1.2}) can be explained as (\ref{1.1}).
In fact, (\ref{1.1}) is equivalent to
\begin{equation}\label{js}
\hat{u}(\xi)=(|x|^{\alpha-n}*[u(1-|u|^2)])^\wedge(\xi)
=C|\xi|^{-\alpha}[u(1-|u|^2)]^\wedge(\xi),
\end{equation}
where $C$ is a positive constant.
By the property of the Riesz potential, we have
\begin{equation}\label{ly}
[(-\Delta)^{\alpha/2}u]^\wedge(\xi)
=C|\xi|^\alpha \hat{u}(\xi),
\end{equation}
where $C$ is another positive constant.
Therefore, the above equality (\ref{js}) amounts to (\ref{1.2}). In addition,
let $u \in E$ be a solution of (\ref{1.2}) with $0<\alpha<2$.
From (\ref{qz}), it follows that (\ref{ly}) is still true.
If the Fourier inversion formula of (\ref{js}) holds,
then $u$ also solves (\ref{1.1}) (if we omit the constants).

\paragraph{Remark 4.}
If $u$ is a finite energy solution of (\ref{1.5}), then \cite{BMR}
shows that
\begin{equation}\label{chou}
\int_{\mathbb{R}^n}|u|^2(1-|u|^2)dx<\infty.
\end{equation}
Therefore, we sometimes call $u$ a finite energy solution
of (\ref{1.1}) if $u$ satisfies (\ref{chou}). Moreover, if
$u$ is uniformly continuous, we can see that either $u \in L^2(\mathbb{R}^n)$
or $1-|u|^2 \in L^1(\mathbb{R}^n)$ by the same argument of (3.9) and (3.10) in \cite{BMR}.
Therefore, if a bounded, uniformly continuous, differentiable function $u$ is
a finite energy solution of (\ref{1.1}), then either $u \equiv 0$, or
$|u(x)| \to 1$ when $|x| \to \infty$.

\paragraph{Acknowledgements.}
This research was supported by NNSF (11871278) of China
and NSF of Jiangsu Education Commission (19KJB110016).


\end{document}